\documentclass[axioms,article,moreauthors,pdftex]{my-mdpi} 

\firstpage{1}
\pubyear{2020}

\history{Submitted: September 4, 2020; Revised: October 4, 2020; Accepted: October 6, 2020}

\pdfoutput=1

\usepackage{epsfig,multicol,graphicx,epstopdf}
\usepackage{mathptmx}

\allowdisplaybreaks 

\Title{Application of Bernoulli Polynomials for Solving Variable-Order Fractional Optimal Control-Affine Problems}

\Author{Somayeh Nemati $^{1}$\orcidA{},
Delfim F. M. Torres $^{2}$*\orcidB{}}

\AuthorNames{Somayeh Nemati and Delfim F. M. Torres}

\address{%
$^{1}$ \quad Department of Mathematics,
Faculty of Mathematical Sciences,
University of Mazandaran, Babolsar, Iran;
s.nemati@umz.ac.ir\\
$^{2}$ \quad Center for Research and Development in Mathematics and Applications (CIDMA),
Department of Mathematics, University of Aveiro, 3810-193 Aveiro, Portugal; delfim@ua.pt}

\corres{Correspondence: delfim@ua.pt; Tel.: +351-234-370-668}

\abstract{We propose two efficient numerical approaches for solving variable-order
fractional optimal control-affine problems. The variable-order fractional
derivative is considered in the Caputo sense, which together with the
Riemann--Liouville integral operator is used in our new techniques.
An accurate operational matrix of variable-order fractional integration
for Bernoulli polynomials is introduced. Our methods proceed as follows.
First, a specific approximation of the differentiation order of the
state function is considered, in terms of Bernoulli polynomials.
Such approximation, together with the initial conditions,
help us to obtain some approximations for the other existing
functions in the dynamical control-affine system. Using these
approximations, and the Gauss--Legendre integration formula,
the problem is reduced to a system of nonlinear algebraic equations.
Some error bounds are then given for the approximate optimal state
and control functions, which allow us to obtain an error bound
for the approximate value of the performance index. We end by
solving some test problems, which demonstrate the high accuracy
of our results.}



\begin{document}


\section{Introduction}
\label{sec:1}

The Bernoulli polynomials, named after Jacob Bernoulli (1654--1705),
occur in the study of many special functions and, in particular,
in relation with fractional calculus, which is a classical area
of mathematical analysis whose foundations were laid by Liouville
in a paper from 1832 and that is nowadays a very active research area
\cite{MR2736622}. One can say that Bernoulli polynomials are a
powerful mathematical tool in dealing with various problems of
dynamical nature \cite{Keshavarz,Bhrawy,Tohidi,Toutounian,Bazm}.
Recently, an approximate method, based on orthonormal Bernoulli's polynomials,
has been developed for solving fractional order differential equations
of Lane--Emden type \cite{MR3957059}, while in \cite{MR3911141}
Bernoulli polynomials are used to numerical solve Fredholm fractional
integro-differential equations with right-sided Caputo derivatives.
Here we are interested in the use of Bernoulli polynomials with respect
to fractional optimal control problems.

An optimal control problem refers to the minimization of a functional
on a set of control and state variables (the performance index)
subject to dynamic constraints on the states and controls.
When such dynamic constraints are described by fractional
differential equations, then one speaks of
fractional optimal control problems (FOCPs) \cite{MR3872489}.
The mathematical theory of fractional optimal control has born
in 1996/97 from practical problems of mechanics and began to be
developed in the context of the fractional calculus of variations
\cite{MR2984893,MR3331286,MR3443073}. Soon after,
fractional optimal control theory became a mature research area,
supported with many applications in engineering and physics. For
the state of the art, see \cite{MR3904404,MR3955185,MR3854267}
and references therein. Regarding the use of Bernoulli polynomials
to numerically solve FOCPs, we refer to \cite{Keshavarz},
where the operational matrices of fractional Riemann--Liouville integration
for Bernoulli polynomials are derived and the properties
of Bernoulli polynomials are utilized, together with
Newton's iterative method, to find approximate solutions to FOCPs.
The usefulness of Bernoulli polynomials for mixed integer-fractional
optimal control problems is shown in \cite{MR3736762}, while the
practical relevance of the methods in engineering
is illustrated in \cite{MR3800767}. Recently, such results
have been generalized for two-dimensional fractional optimal control problems,
where the control system is not a fractional ordinary differential equation
but a fractional partial differential equation \cite{MR3927863}.
Here we are the first to develop a numerical method,
based on Bernoulli polynomials, for FOCPs of variable-order.

The variable-order fractional calculus was introduced in 1993
by Samko and Ross and deals with operators of order $\alpha$,
where $\alpha$ is not necessarily a constant but a function
$\alpha(t)$ of time \cite{MR1421643}. With this extension,
numerous applications have been found in physics, mechanics,
control, and signal processing \cite{MR1926468,MR3930479,MR3958579,MR3060420,MR3955314}.
For the state of the art on variable-order fractional optimal control
we refer the interested reader to the book \cite{Almeida} and the articles
\cite{MR3864326,MR3907973}. To the best of our knowledge, numerical methods
based on Bernoulli polynomials for such kind of FOCPs are not available
in the literature. For this reason, in this work we focus
on the following variable-order fractional
optimal control-affine problem (FOC-AP):
\begin{equation}
\label{1.1}
\min~J=\int_0^1\phi(t,x(t),u(t))dt
\end{equation}
subject to the control-affine dynamical system
\begin{equation}
\label{1.2}
{_0^CD}_t^{\alpha(t)}x(t)=\varphi\left(t,x(t),{_0^CD}_t^{\alpha_1(t)}x(t),
\dots,{_0^CD}_t^{\alpha_s(t)}x(t)\right)+b(t)u(t)
\end{equation}
and the initial conditions
\begin{equation}
\label{1.3}
x^{(i)}(0)=x_0^i,\quad i=0,1,\ldots,n,
\end{equation}
where $\phi$ and $\varphi$ are smooth functions of their arguments,
$b\neq 0$, $n$ is a positive integer number such that for all $t\in[0,1]$,
$0<\alpha_1(t)<\alpha_2(t)<\ldots <\alpha_s(t)<\alpha(t)\leq n$, and
${_0^CD_t^{\alpha(t)}}$ is the (left) fractional derivative of variable-order
defined in the Caputo sense. We employ two spectral methods based on Bernoulli
polynomials in order to obtain numerical solutions to problem \eqref{1.1}--\eqref{1.3}.
Our main idea consists of reducing the problem to a system of nonlinear algebraic equations.
To do this, we introduce an accurate operational matrix of variable-order fractional integration,
having Bernoulli polynomials as basis vectors.

The paper is organized as follows. In Section~\ref{sec:2}, the variable-order fractional
calculus is briefly reviewed and some properties of the Bernoulli polynomials are recalled.
A new operational matrix of variable-order is introduced for the Bernoulli basis functions
in Section~\ref{sec:03}. Section~\ref{sec:3} is devoted to two new numerical approaches
based on Bernoulli polynomials for solving problem \eqref{1.1}--\eqref{1.3}. In Section~\ref{sec:4},
some error bounds are proved. Then, in Section~\ref{sec:5}, some FOC-APs are solved
using the proposed methods. Finally, concluding remarks are given in Section~\ref{sec:6}.


\section{Preliminaries}
\label{sec:2}

In this section, a brief review on necessary definitions and properties of the
variable-order fractional calculus is presented. Furthermore, Bernoulli
polynomials and some of their properties are recalled.


\subsection{The variable-order fractional calculus}

The two most commonly used definitions in fractional calculus are the
Riemann--Liouville integral and the Caputo derivative. Here, we deal with
generalizations of these two definitions, which allow the order of the
fractional operators to be of variable-order.

\begin{Definition}[See, e.g., \cite{Almeida}]
The left Riemann--Liouville fractional integral of order $\alpha(t)$ is defined by
\begin{equation*}
_0I_t^{\alpha(t)}y(t)=\frac{1}{\Gamma{(\alpha(t))}}
\int_0^t(t-s)^{\alpha(t)-1}y(s)ds,\quad t>0,
\end{equation*}
where $\Gamma$ is the Euler gamma function, that is,
\begin{equation*}
\Gamma(t)=\int_0^{\infty}\tau^{t-1}\exp(-\tau) d\tau,
\quad t>0.
\end{equation*}
\end{Definition}

\begin{Definition}[See, e.g., \cite{Almeida}]
The left Caputo fractional derivative
of order $\alpha(t)$ is defined by
\begin{align*}
&{_0^CD}_t^{\alpha(t)}y(t)=\frac{1}{\Gamma(n-\alpha(t))}
\int_0^t(t-s)^{n-\alpha(t)-1}y^{(n)}(s)ds,\quad n-1<\alpha(t)<n,\\
&{_0^CD}_t^{\alpha(t)}y(t)=y^{(n)}(t),\quad \alpha(t)=n.
\end{align*}
\end{Definition}

For $0\leq \beta(t)<\alpha(t)\leq n$, $n\in\mathbb{N}$, $\gamma>0$, and $\nu>-1$,
some useful properties of the Caputo derivative and Riemann--Liouville
fractional integral are as follows \cite{Almeida}:
\begin{equation}
\label{2.1}
_0I_t^{\alpha(t)}t^\nu
=\frac{\Gamma(\nu+1)}{\Gamma{(\nu+1+\alpha(t))}}t^{\nu+\alpha(t)},
\end{equation}
\begin{equation}
\label{2.2}
_0I_t^{\gamma}({_0^CD}_t^{\gamma}y(t))=y(t)
-\sum_{i=0}^{\lceil\gamma\rceil-1}y^{(i)}(0)\frac{t^i}{i!},
\quad t>0,
\end{equation}
\begin{equation}
\label{2.3}
_0I_t^{n-\alpha(t)}(y^{(n)}(t))={_0^CD}_t^{\alpha(t)}y(t)
-\sum_{i=\lceil\alpha(t)\rceil}^{n-1}y^{(i)}(0)
\frac{t^{i-\alpha(t)}}{\Gamma(i+1-\alpha(t))},
\quad t>0,
\end{equation}
\begin{equation}
\label{2.4}
_0I_t^{\alpha(t)-\beta(t)}({_0^CD}_t^{\alpha(t)}y(t))
={_0^CD}_t^{\beta(t)}y(t)-\sum_{i=\lceil\beta(t)\rceil}^{\lceil
\alpha(t)\rceil-1}y^{(i)}(0)\frac{t^{i-\beta(t)}}{\Gamma(i+1-\beta(t))},
\quad t>0,
\end{equation}
where $\lceil \cdot \rceil$ is the ceiling function.


\subsection{Bernoulli polynomials}

The set of Bernoulli polynomials, $\{\beta_m(t)\}_{m=0}^{\infty}$, consists
of a family of independent functions that builds a complete basis for the space
$L^2[0,1]$ of all square integrable functions on the interval $[0,1]$.
These polynomials are defined as
\begin{equation}
\label{2.5}
\beta_m(t)=\sum_{i=0}^m \binom{m}{i}b_{m-i}t^i,
\end{equation}
where $b_k$, $k=0,1,\ldots,m$, are the Bernoulli numbers \cite{Costabile}.
These numbers are seen in the series expansion of trigonometric functions
and can be given by the following identity \cite{Arfken}:
\begin{equation*}
\frac{t}{e^t-1}
=\sum_{i=0}^\infty b_i \frac{t^i}{i!}.
\end{equation*}
Thus, the first few Bernoulli numbers are given by
\begin{eqnarray*}
b_0=1,\quad
b_1=-\frac{1}{2},\quad
b_2=\frac{1}{6},\quad
b_3 = 0,\quad
b_4=-\frac{1}{30},\quad
b_5 = 0,\quad
b_6=\frac{1}{42}.
\end{eqnarray*}
Furthermore, the first five Bernoulli polynomials are
\begin{align*}
&\beta_0(t)=1,\\
&\beta_1(t)=t-\frac{1}{2},\\
&\beta_2(t)=t^2-t+\frac{1}{6},\\
&\beta_3(t)=t^3-\frac{3}{2}t^2+\frac{1}{2}t,\\
&\beta_4(t)=t^4-2t^3+t^2-\frac{1}{30}.
\end{align*}

For an arbitrary function $x\in L^2[0,1]$, we can write
\begin{equation*}
x(t) = \sum_{m=0}^\infty a_m\beta_m(t).
\end{equation*}
Therefore, an approximation of the function $x$ can be given by
\begin{equation}
\label{2.7}
x(t)\simeq x_M(t) = \sum_{m=0}^M a_m\beta_m(t)=A^TB(t),
\end{equation}
where
\begin{equation}\label{2.8}
B(t)=[\beta_0(t),\beta_1(t),\ldots,\beta_M(t)]^T
\end{equation}
and
\begin{equation*}
A=[a_0,a_1,\ldots,a_M]^T.
\end{equation*}
The vector $A$ in \eqref{2.7} is called the coefficient vector
and can be calculated by the formula (see \cite{Keshavarz})
\begin{equation*}
A=D^{-1}\langle x(t),B(t)\rangle,
\end{equation*}
where $\langle \cdot,\cdot \rangle$ is the inner product,
defined for two arbitrary functions $f,g\in L^2 [0,1]$ as
\begin{equation*}
\langle f(t),g(t) \rangle=\int_0^1f(t)g(t)dt,
\end{equation*}
and $D=\langle B(t),B(t)\rangle$ is calculated using the following
property of Bernoulli polynomials \cite{Arfken}:
\begin{equation*}
\int_0^1\beta_i(t)\beta_j(t)dt=(-1)^{i-1}\frac{i!j!}{(i+j)!}b_{i+j},
\quad i,j\geq 1.
\end{equation*}
It should be noted that
\begin{equation*}
X=span\left\{\beta_0(t),\beta_1(t),\ldots,\beta_M(t)\right\}
\end{equation*}
is a finite dimensional subspace of $L^2[0,1]$ and $x_M$,
given by \eqref{2.7}, is the best approximation of function $x$ in $X$.


\section{Operational matrix of variable-order fractional integration}
\label{sec:03}

In this section, we introduce an accurate operational matrix
of variable-order fractional integration for Bernoulli functions.
To this aim, we rewrite the Bernoulli basis vector $B$ given
by \eqref{2.8} in terms of the Taylor basis functions as follows:
\begin{equation}
\label{3.1}
B(t)=Q\mathbb{T}(t),
\end{equation}
where $\mathbb{T}$ is the Taylor basis vector given by
\begin{equation*}
\mathbb{T}(t)=\left[1,t,t^2,\ldots,t^M\right]^T
\end{equation*}
and $Q$ is the change-of-basis matrix, which is obtained using \eqref{2.5} as
\begin{equation*}
Q=\left[
\begin{array}{ccccccc}
1& 0 & 0 & 0 & 0 &\ldots & 0\\
-\frac{1}{2} & 1 & 0 & 0 & 0&\ldots & 0\\
\frac{1}{6} & -1 & 1 & 0 & 0&\ldots & 0\\
0 &\frac{1}{2}& -\frac{3}{2} & 1 &0 &\ldots & 0\\
\vdots & \vdots& \vdots & \vdots &\vdots& &  \vdots\\
b_M & \binom{M}{1}b_{M-1}  & \binom{M}{2}b_{M-2}
& \binom{M}{3}b_{M-3} & \binom{M}{4}b_{M-4}&\ldots & 1\\
\end{array}
\right].
\end{equation*}
Since $Q$ is nonsingular, we can write
\begin{equation}\label{3.2}
\mathbb{T}(t)=Q^{-1}B(t).
\end{equation}
By considering \eqref{3.1} and applying
the left Riemann--Liouville fractional integral
operator of order $\alpha(t)$ to the vector $B(t)$, we get that
\begin{equation}
\label{3.3}
_0I_t^{\alpha(t)} B(t)={_0I_t^{\alpha(t)}}(Q\mathbb{T}(t))
=Q ({_0I_t^{\alpha(t)}}\mathbb{T}(t))=QS_t^{\alpha(t)}\mathbb{T}(t),
\end{equation}
where $S_t^{\alpha(t)}$ is a diagonal matrix,
which is obtained using \eqref{2.1} as follows:
\begin{equation*}
S_t^{\alpha(t)}=\left[
\begin{array}{cccccc}
\frac{1}{\Gamma(1+\alpha(t))}t^{\alpha(t)} & 0 & 0 &0&\cdots &0\\
0 & \frac{1}{\Gamma(2+\alpha(t))}t^{\alpha(t)} & 0 &0&\cdots &0\\
0 & 0 & \frac{2}{\Gamma(3+\alpha(t))}t^{\alpha(t)} &0&\cdots &0\\
\vdots & \vdots & \vdots & \vdots& & \vdots\\
0 & 0 & 0 & 0 & \cdots & \frac{\Gamma(M+1)}{\Gamma{(M+1+\alpha(t))}}t^{\alpha(t)}\\
\end{array}
\right].
\end{equation*}
Finally, by using \eqref{3.2} in \eqref{3.3}, we have
\begin{equation}
\label{3.4}
_0I_t^{\alpha(t)} B(t)=QS_t^{\alpha(t)}Q^{-1}B(t)=P_t^{\alpha(t)}B(t),
\end{equation}
where $P_t^{\alpha(t)}=QS_t^{\alpha(t)}Q^{-1}$ is a matrix of dimension $(M+1)\times (M+1)$,
which we call the operational matrix of variable-order fractional integration $\alpha(t)$
for Bernoulli functions. Since $Q$ and $Q^{-1}$ are lower triangular matrices and $S_t^{\alpha(t)}$
is a diagonal matrix, $P_t^{\alpha(t)}$ is also a lower triangular matrix.
In the particular case of $M=2$, one has
\begin{equation*}
P_t^{\alpha(t)}=\left[
\begin{array}{ccc}
\frac{1}{\Gamma (\alpha (t)+1)}t^{\alpha (t)} & 0 & 0 \\
\left(\frac{1}{2 \Gamma (\alpha (t)+2)}
-\frac{1}{2 \Gamma (\alpha (t)+1)}\right) t^{\alpha (t)}
& \frac{1}{\Gamma (\alpha (t)+2)}t^{\alpha (t)} & 0 \\
\left(\frac{1}{6 \Gamma (\alpha (t)+1)}-\frac{1}{2 \Gamma (\alpha (t)+2)}
+\frac{2 }{3 \Gamma (\alpha (t)+3)}\right)t^{\alpha (t)}
& \left(\frac{2}{\Gamma (\alpha (t)+3)}
-\frac{1}{\Gamma (\alpha (t)+2)} \right)t^{\alpha (t)}
& \frac{2}{\Gamma (\alpha (t)+3)} t^{\alpha (t)}
\end{array}
\right].
\end{equation*}


\section{Methods of solution}
\label{sec:3}

In this section, we propose two approaches for solving problem
\eqref{1.1}--\eqref{1.3}. To do this, first we introduce
\begin{equation*}
n=\max_{0< t\leq 1}\left\{\lceil \alpha(t)\rceil\right\}.
\end{equation*}
Then, we may use the following two approaches to find approximations
for the state and control functions, which optimize the performance index.


\subsection{Approach I}

In our first approach, we consider an approximation of the $n$th order derivative
of the unknown state function $x$ using Bernoulli polynomials. Set
\begin{equation}
\label{4.1}
x^{(n)}(t)=A^TB(t),
\end{equation}
where $A$ is a $1\times (M+1)$ vector with unknown elements
and $B$ is the Bernoulli basis vector given by \eqref{2.8}. Then,
using the initial conditions given in \eqref{1.3},
and equations \eqref{2.2}, \eqref{3.4}, and \eqref{4.1}, we get
\begin{equation}
\label{4.2}
\begin{split}
x(t)&={_0I_t^{n}}(x^{(n)}(t))+\sum_{i=0}^{n-1}x^{(i)}(0)\frac{t^i}{i!}\\
&=A^T({_0I_t^{n}}B(t))+\sum_{i=0}^{n-1}x_0^{i}\frac{t^i}{i!}\\
&=A^TP_t^{n}B(t)+\sum_{i=0}^{n-1}x_0^{i}\frac{t^i}{i!}.
\end{split}
\end{equation}
Moreover, using \eqref{2.3}, \eqref{3.4}, and \eqref{4.1}, we get
\begin{equation}
\label{4.3}
{_0^CD}_t^{\alpha(t)}x(t)=A^TP_t^{n-\alpha(t)}B(t)
+\sum_{i=\lceil\alpha(t)\rceil}^{n-1}
x_0^{i}\frac{t^{i-\alpha(t)}}{\Gamma(i+1-\alpha(t))}:=F[A,t]
\end{equation}
and
\begin{equation}
\label{4.4}
{_0^CD}_t^{\alpha_j(t)}x(t)=A^TP_t^{n-\alpha_j(t)}B(t)
+\sum_{i=\lceil\alpha_j(t)\rceil}^{n-1}x_0^{i}
\frac{t^{i-\alpha_j(t)}}{\Gamma(i+1-\alpha_j(t))}:=F_j[A,t],
\quad j=1,\ldots,s.
\end{equation}
By substituting \eqref{4.2}--\eqref{4.4} into the control-affine dynamical
system given by \eqref{1.2}, we obtain an approximation
of the control function as follows:
\begin{equation}
\label{4.5}
u(t)=\frac{1}{b(t)}\left[F[A,t]-\varphi\left(t,A^TP_t^{n}B(t)
+\sum_{i=0}^{n-1}x_0^{i}\frac{t^i}{i!},F_1[A,t],\dots,F_s[A,t]\right)\right].
\end{equation}
Taking into consideration \eqref{4.2} and \eqref{4.5} in the performance index $J$, we have
\begin{equation*}
J[A]=\int_0^1{\phi\left(t,A^TP_t^{n}B(t)+\sum_{i=0}^{n-1}x_0^{i}\frac{t^i}{i!},
\frac{1}{b(t)}\left[F[A,t]-\varphi\left(t,A^TP_t^{n}B(t)+\sum_{i=0}^{n-1}x_0^{i}
\frac{t^i}{i!},F_1[A,t],\dots,F_s[A,t]\right)\right]\right)dt}.
\end{equation*}
For the sake of simplicity, we introduce
\begin{equation*}
G[A,t]=\phi\left(t,A^TP_t^{n}B(t)+\sum_{i=0}^{n-1}x_0^{i}\frac{t^i}{i!},
\frac{1}{b(t)}\left[F[A,t]-\varphi\left(t,A^TP_t^{n}B(t)
+\sum_{i=0}^{n-1}x_0^{i}\frac{t^i}{i!},F_1[A,t],\dots,F_s[A,t]\right)\right]\right).
\end{equation*}
In many applications, it is difficult to compute the integral of function $G[A,t]$.
Therefore, it is recommended to use a suitable numerical integration formula. Here,
we use the Gauss--Legendre quadrature formula to obtain
\begin{equation}
\label{index}
J[A]\simeq \frac{1}{2}\sum_{i=1}^N\omega_iG\left[A,\frac{t_i+1}{2}\right],
\end{equation}
where $t_i$, $i=1,2,\ldots,N$, are the zeros of the Legendre polynomial of degree $N$,
$P_N(t)$, and $\omega_i$ are the corresponding weights \cite{Shen1}, which are given by
\begin{equation}
\label{wei}
\omega_i=\frac{2}{\left( \frac{d}{dt}P_N(t_i)\right)^2(1-t_i^2)},\quad i=1,\ldots,N.
\end{equation}
Finally, the first order necessary condition for the optimality of the performance index implies
\begin{equation*}
\frac{\partial J[A]}{\partial A}=0,
\end{equation*}
which gives a system of $M+1$ nonlinear algebraic equations in terms of the $M+1$
unknown elements of the vector $A$. By solving this system, approximations of the
optimal state and control functions are, respectively, given by \eqref{4.2} and \eqref{4.5}.


\subsection{Approach II}

In our second approach, we set
\begin{equation}
\label{4.6}
{_0^CD}_t^{\alpha(t)}x(t)=A^TB(t).
\end{equation}
Then, using \eqref{2.4} with $\beta(t)\equiv 0$, we obtain that
\begin{equation}
\label{4.7}
\begin{split}
x(t)&={_0I_t^{\alpha(t)}}({_0^CD}_t^{\alpha(t)}x(t))
+\sum_{i=0}^{\lceil\alpha(t)\rceil-1}x^{(i)}(0)\frac{t^{i}}{\Gamma(i+1)}\\
&=A^T({_0I_t^{\alpha(t)}}B(t))+\sum_{i=0}^{\lceil\alpha(t)\rceil-1}x_0^{i}\frac{t^i}{i!}\\
&=A^TP_t^{\alpha(t)}B(t)+\sum_{i=0}^{\lceil\alpha(t)\rceil-1}x_0^{i}\frac{t^i}{i!}.
\end{split}
\end{equation}
Furthermore, we get
\begin{equation}
\label{4.8}
{_0^CD}_t^{\alpha_j(t)}x(t)=A^TP_t^{\alpha(t)-\alpha_j(t)}B(t)
+\sum_{i=\lceil\alpha_j(t)\rceil}^{\lceil\alpha(t)\rceil-1}x_0^{i}
\frac{t^{i-\alpha_j(t)}}{\Gamma(i+1-\alpha_j(t))}:=F_j[A,t],
\quad j=1,\ldots,s.
\end{equation}
Taking \eqref{4.6}--\eqref{4.8} into consideration, equation \eqref{1.2} gives
\begin{equation}
\label{4.9}
u(t)=\frac{1}{b(t)}\left[A^TB(t)-\varphi\left(t,A^TP_t^{\alpha(t)}B(t)
+\sum_{i=0}^{\lceil\alpha(t)\rceil-1}x_0^{i}\frac{t^i}{i!},F_1[A,t],
\dots,F_s[A,t]\right)\right].
\end{equation}
By substituting the approximations given by \eqref{4.7} and \eqref{4.9}
into the performance index, we get
\begin{multline*}
J[A]=\int_0^1\phi\left(t,A^TP_t^{\alpha(t)}B(t)
+\sum_{i=0}^{\lceil\alpha(t)\rceil-1}x_0^{i}\frac{t^i}{i!},\right.\\
\left.\frac{1}{b(t)}\left[A^TB(t)-\varphi\left(t,A^TP_t^{\alpha(t)}B(t)
+\sum_{i=0}^{\lceil\alpha(t)\rceil-1}x_0^{i}
\frac{t^i}{i!},F_1[A,t],\dots,F_s[A,t]\right)\right]\right)dt.
\end{multline*}
By introducing
\begin{multline*}
G[A,t]=\phi\left(t,A^TP_t^{\alpha(t)}B(t)+\sum_{i=0}^{\lceil\alpha(t)\rceil-1}
x_0^{i}\frac{t^i}{i!},\right.\\
\left.\frac{1}{b(t)}\left[A^TB(t)-\varphi\left(t,A^TP_t^{\alpha(t)}B(t)
+\sum_{i=0}^{\lceil\alpha(t)\rceil-1}x_0^{i}\frac{t^i}{i!},
F_1[A,t],\dots,F_s[A,t]\right)\right]\right),
\end{multline*}
then this approach continues in the same way of finding the unknown parameters
of the vector $A$ as in Approach~I.


\section{Error bounds}
\label{sec:4}

The aim of this section is to give some error bounds for the numerical solution obtained
by the proposed methods of Section~\ref{sec:3}. We present the error discussion for
Approach~II, which can then be easily extended to Approach~I.

Suppose that $x^*$ is the optimal state function of problem \eqref{1.1}--\eqref{1.3}.
Let $f(t):={_0^CD}_t^{\alpha(t)}x^*(t)$ with $f(t)\in H^\mu (0,1)$
($H^\mu (0,1)$ is a Sobolev space \cite{Canuto}). According to our numerical method,
$f_M(t)=A^TB(t)$ is the best approximation of function $f$ in terms of the Bernoulli polynomials,
that is,
\begin{equation*}
\forall g\in X,~\|f-f_M\|_2\leq \|f-g\|_2.
\end{equation*}
We recall the following lemma from \cite{Canuto}.

\begin{Lemma}[See \cite{Canuto}]
\label{lem1}
Assume that $f\in H^\mu (0,1)$ with $\mu\geq 0$. Let $L_M(f)\in X$
be the truncated shifted Legendre series of $f$. Then,
\begin{equation*}
\|f-L_M(f)\|_2\leq CM^{-\mu}|f|_{H^{\mu;M}(0,1)},
\end{equation*}
where
\begin{equation*}
|f|_{H^{\mu;M}(0,1)}=\left(\sum_{j=\min\{\mu,M+1\}}^\mu
\|f^{(j)}\|_2^2\right)^{\frac{1}{2}}
\end{equation*}
and $C$ is a positive constant independent of function $f$ and integer $M$.
\end{Lemma}

Since the best approximation of function $f$ in the subspace $X$ is unique
and $f_M$ and $L_M(f)$ are both the best approximations of $f$ in $X$,
we have $f_M=L_M(f)$. Therefore, we get that
\begin{equation}
\label{bound}
\|f-f_M\|_2\leq CM^{-\mu}|f|_{H^{\mu;M}(0,1)}.
\end{equation}
Hereafter, $C$ denotes a positive constant independent of $M$ and $n$.

\begin{Theorem}
\label{th1}
Suppose $x^*$ to be the exact optimal state function of problem \eqref{1.1}--\eqref{1.3}
such that $f(t):={_0^CD}_t^{\alpha(t)}x^*(t)\in H^\mu (0,1)$, with $\mu\geq 0$,
and $\tilde{x}$ be its approximation given by \eqref{4.7}. Then,
\begin{equation}
\label{bound1}
\|x^*(t)-\tilde{x}(t)\|_2
\leq CM^{-\mu}|f|_{H^{\mu;M}(0,1)}.
\end{equation}
\end{Theorem}

\begin{proof}
Let $Y=L^2[0,1]$ and $_0I_t^{\alpha(t)}:Y\rightarrow Y$ be the variable-order
Riemann--Liouville integral operator. By definition of the norm for operators,
we have
\begin{equation*}
\|_0I_t^{\alpha(t)}\|_2=\sup_{\|g\|_2=1}\|_0I_t^{\alpha(t)}g\|_2.
\end{equation*}
In order to prove the theorem, first we show that the operator $_0I_t^{\alpha(t)}$ is bounded.
Since $\|g\|_2=1$, using Schwarz's inequality, we get
\begin{align*}
\left\|_0I_t^{\alpha(t)}g\right\|_2
&=\left\|\frac{1}{\Gamma(\alpha(t))}\int_0^t(t-s)^{\alpha(t)-1}g(s)ds\right\|_2\\
&\leq \|g\|_2\left\|\frac{1}{\Gamma(\alpha(t))}\int_0^t(t-s)^{\alpha(t)-1}ds\right\|_2\\
&=\left\|\frac{t^{\alpha(t)}}{\Gamma(\alpha(t)+1)}\right\|_2\\
&\leq C,
\end{align*}
where we have used the assumption $\alpha(t)>0$, which gives $t^{\alpha(t)}<1$
for $0<t\leq 1$, and a particular property of the Gamma function, which is $\Gamma(t)>0.8$.
Therefore, $_0I_t^{\alpha(t)}$ is bounded. Now, using \eqref{bound}, and taking into account
\eqref{2.4} and \eqref{4.7}, we obtain that
\begin{align*}
\left\|x^*(t)-\tilde{x}(t)\right\|_2
&=\left\|{_0I_t^{\alpha(t)}}f(t)
+\sum_{i=0}^{\lceil\alpha(t)\rceil-1}x^{(i)}(0)\frac{t^{i}}{\Gamma(i+1)}
-\left( {_0I_t^{\alpha(t)}}(A^TB(t))+\sum_{i=0}^{\lceil\alpha(t)\rceil-1}
x_0^{i}\frac{t^i}{i!}\right)\right\|_2\\
&=\left\|{_0I_t^{\alpha(t)}}(f(t)-A^TB(t))\right\|_2\\
&~\leq \left\|{_0I_t^{\alpha(t)}} \right\|_2\left\|f(t)-A^TB(t)\right\|_2\\
&\leq C M^{-\mu}|f|_{H^{\mu;M}(0,1)}.
\end{align*}
The proof is complete.
\end{proof}

\begin{Remark}
Since we have $\alpha(t)-\alpha_j(t)>0$, $j=1,2,\ldots,s$,
with a similar argument it can be shown that
\begin{equation*}
\left\|{_0^CD}_t^{\alpha_j(t)}x^*(t)-\left(A^TP_t^{\alpha(t)
-\alpha_j(t)}B(t)+\sum_{i=\lceil\alpha_j(t)\rceil}^{\lceil\alpha(t)\rceil-1}
x_0^{i}\frac{t^{i-\alpha_j(t)}}{\Gamma(i+1-\alpha_j(t))}\right)\right\|_2
\leq C M^{-\mu}|f|_{H^{\mu;M}(0,1)}.
\end{equation*}
\end{Remark}

With the help of Theorem~\ref{th1}, we obtain the following result
for the error of the optimal control function. For simplicity, suppose
that in the control-affine dynamical system given by \eqref{1.2}
the function $\varphi$ appears as $\varphi:=\varphi(t,x)$
(cf. Remark~\ref{rem1}).

\begin{Theorem}
\label{th2}
Suppose that the assumptions of Theorem~\ref{th1} are fulfilled.
Let $u^*$ and $\tilde{u}$ be the exact and approximate optimal control functions,
respectively. If $\varphi:\mathbb{R}^2\longrightarrow\mathbb{R}$ satisfies
a Lipschitz condition with respect to the second argument, then
\begin{equation}
\label{bound2}
\|u^*(t)-\tilde{u}(t)\|_2\leq C M^{-\mu}|f|_{H^{\mu;M}(0,1)}.
\end{equation}
\end{Theorem}

\begin{proof}
Using equation \eqref{1.2}, the exact optimal control function is given by
\begin{equation}
\label{5.1}
u^*(t)=\frac{1}{b(t)}\left(f(t)-\varphi\left(t,x^*(t)\right)\right)
\end{equation}
and the approximate control function obtained by our Approach~II is given by
\begin{equation}
\label{5.2}
\tilde{u}(t)=\frac{1}{b(t)}\left(A^TB(t)-\varphi\left(t,\tilde{x}(t)\right)\right).
\end{equation}
By subtracting \eqref{5.2} from \eqref{5.1}, we get
\begin{equation}
\label{5.3}
u^*(t)-\tilde{u}(t)=\frac{1}{b(t)}\left(f(t)-\varphi\left(t,x^*(t)\right)
-A^TB(t)+\varphi\left(t,\tilde{x}(t)\right)\right).
\end{equation}
Since the function $\varphi$ satisfies a Lipschitz condition with respect
to the second variable, there exists a positive constant $K$ such that
\begin{equation*}
|\varphi(t,{\bf x}_1)-\varphi(t,{\bf x}_2)|<K|{\bf x}_1-{\bf x}_2|.
\end{equation*}
Therefore, using \eqref{bound} and \eqref{bound1}, and also
taking into account $b(t)\neq 0$, we have
\begin{equation*}
\|u^*(t)-\tilde{u}(t)\|_2\leq \frac{1}{\|b(t)\|_2}\left(\left\|f(t)
-A^TB(t)\right\|_2+K\left\|x^*(t)-\tilde{x}(t)\right\|_2\right)
\leq C M^{-\mu}|f|_{H^{\mu;M}(0,1)},
\end{equation*}
which yields \eqref{bound2}.
\end{proof}

\begin{Remark}
\label{rem1}
For the general case $\varphi:=\varphi(t,x,x_1,\ldots,x_s)$, the same result
of Theorem~\ref{th2} can be easily obtained by assuming that $\varphi$
satisfies Lipschitz conditions with respect to the variables $x$, $x_1$, \ldots, $x_s$.
\end{Remark}

As a result of Theorems~\ref{th1} and \ref{th2}, we obtain an error bound for
the approximate value of the optimal performance index $J$ given by \eqref{index}.
First, we recall the following lemma in order to obtain the error
of the Gauss--Legendre quadrature rule.

\begin{Lemma}[See \cite{Shen1}]
Let $g$ be a given sufficiently smooth function. Then,
the Gauss--Legendre quadrature rule is given by
\begin{equation}
\label{int}
\int_{-1}^1g(t)dt=\sum_{i=1}^N \omega_i g(t_i)+E_N(g),
\end{equation}
where $t_i$, $i=1,\ldots,N$, are the roots of the Legendre polynomial of degree $N$,
and $\omega_i$ are the corresponding weights given by \eqref{wei}.
In \eqref{int}, $E_N(g)$ is the error term, which is given as follows:
\begin{equation*}
E_N(g)=\frac{2^{2N+1}(N!)^4}{(2N+1)[(2N!)]^3}g^{2N}(\eta),\quad \eta\in(-1,1).
\end{equation*}
\end{Lemma}

Now, by considering the assumptions of Theorems~\ref{th1} and \ref{th2},
we prove the following result.

\begin{Theorem}
\label{th3}
Let $J^*$ be the exact value of the optimal performance index $J$
in problem \eqref{1.1}--\eqref{1.3} and $\tilde{J}$ be its approximation
given by \eqref{index}. Suppose that the function $\phi:\mathbb{R}^3
\longrightarrow \mathbb{R}$ is a sufficiently smooth function with
respect to all its variables and satisfies Lipschitz conditions
with respect to its second and third arguments, that is,
\begin{equation}
\label{L1}
|\phi(t,{\bf x}_1,u)-\phi(t,{\bf x}_2,u)|\leq K_1|{\bf x}_1-{\bf x}_2|
\end{equation}
and
\begin{equation}
\label{L2}
|\phi(t,x,{\bf u}_1)-\phi(t,x,{\bf u}_2)|\leq K_1|{\bf u}_1-{\bf u}_2|,
\end{equation}
where $K_1$ and $K_2$ are real positive constants. Then, there exist
positive constants $C_1$ and $C_2$ such that
\begin{equation}
\label{bound:r}
\left|J^*-\tilde{J}\right|
\leq C_1M^{-\mu}|f|_{H^{\mu;M}(0,1)}
+C_2 \frac{(N!)^4}{(2N+1)[(2N!)]^3}.
\end{equation}
\end{Theorem}

\begin{proof}
Using \eqref{index} and \eqref{int}, we have
\begin{equation}
\label{J}
\tilde{J}=\frac{1}{2}\sum_{i=1}^N\omega_i\phi\left(\frac{t_i+1}{2},
\tilde{x}\left(\frac{t_i+1}{2}\right),\tilde{u}\left(\frac{t_i+1}{2}\right)\right)
=\int_0^1\phi\left(t,\tilde{x}(t),\tilde{u}(t)\right)dt-\xi_N,
\end{equation}
where
\begin{equation*}
\xi_N=\left(\frac{1}{2}\right)\frac{2^{2N+1}(N!)^4}{(2N+1)[(2N!)]^3}\left(
\frac{1}{2} \right)^{2N}\left.\frac{\partial^{2N}\phi\left(t,\tilde{x}(t),
\tilde{u}(t)\right)}{\partial t^{2N}}\right|_{t=\eta}=\frac{(N!)^4}{(2N+1)[(2N!)]^3}\left.
\frac{\partial^{2N}\phi\left(t,\tilde{x}(t),\tilde{u}(t)\right)}{\partial t^{2N}}\right|_{t=\eta}
\end{equation*}
for $\eta\in (0,1)$. Therefore, taking into consideration \eqref{L1}--\eqref{J}, we get
\begin{align*}
\left|J^*-\tilde{J}\right|
&=\left|\int_0^1\phi(t,x^*(t),u^*(t))dt
-\int_0^1\phi\left(t,\tilde{x}(t),\tilde{u}(t)\right)dt+\xi_N\right|\\
&=\left|\int_0^1\phi(t,x^*(t),u^*(t))dt-\int_0^1\phi(t,\tilde{x}(t),
u^*(t))dt+\int_0^1\phi(t,\tilde{x}(t),u^*(t))dt
-\int_0^1\phi\left(t,\tilde{x}(t),\tilde{u}(t)\right)dt+\xi_N\right|\\
&\leq K_1\int_0^1\left|x^*(t)-\tilde{x}(t)\right|dt+K_2
\int_0^1\left|u^*(t)-\tilde{u}(t)\right|dt+\max_{0<t<1}|\xi_N|\\
&\leq C_1M^{-\mu}|f|_{H^{\mu;M}(0,1)}+C_2 \frac{(N!)^4}{(2N+1)[(2N!)]^3},
\end{align*}
where we have used the property of equivalence of $L^1$ and $L^2$-norms and
\begin{equation*}
C_2=\max_{0<t<1}\left|\frac{\partial^{2N}\phi\left(t,\tilde{x}(t),
\tilde{u}(t)\right)}{\partial t^{2N}}\right|.
\end{equation*}
The proof is complete.
\end{proof}

\begin{Remark}
A similar error discussion can be considered for Approach~I 
by setting $f(t):={x^*}^{(n)}(t)$ with $f(t)\in H^\mu(0,1)$ 
and taking into account the fact that the operators $I^n$, 
$I^{\alpha(t)}$ and $I^{\alpha_j(t)}$, for $j=1,2,\ldots,s$, 
are bounded.
\end{Remark}

\begin{Remark}
In practice, since the exact control and state functions that minimize
the performance index are unknown, in order to reach a given specific accuracy
$\epsilon$ for these functions, we increase the number of basis functions
(by increasing $M$) in our implementation, such that
\begin{equation*}
\max_{1\leq i\leq M}\left|F[A,t_i]-\varphi\left(t_i,\tilde{x}(t_i),F_1[A,t_i],
\dots,F_s[A,t_i]\right)-b(t_i)\tilde{u}(t_i)\right|<\epsilon
\quad \text{(Approach~I)},
\end{equation*}
and
\begin{equation*}
\max_{1\leq i\leq M}\left|A^TB(t_i)-\varphi\left(t_i,\tilde{x}(t_i),F_1[A,t_i],
\dots,F_s[A,t_i]\right)-b(t_i)\tilde{u}(t_i)\right|<\epsilon
\quad \text{(Approach~II)},
\end{equation*}
where
\begin{equation*}
t_i=\frac{i}{M+1},\quad i=1,2,\ldots,M.
\end{equation*}
\end{Remark}


\section{Test problems}
\label{sec:5}

In this section, some FOC-APs are included and solved by the proposed methods,
in order to illustrate the accuracy and efficiency of the new techniques.
In our implementation, the method was carried out using
\textsf{Mathematica 12}. Furthermore, we have used $N=14$
in employing the Gauss--Legendre quadrature formula.

\begin{Example}
\label{ex1}
As first example, we consider the following variable-order FOC-AP:
\begin{equation}
\label{6.1}
\min~J=\int_0^1\left[\left(x(t)-t^2\right)^2+\left(u(t)
-\frac{1}{\Gamma(3-\alpha(t))}t^{2-\alpha(t)} e^{-t}
+\frac{1}{2}e^{t^2-t}\right)^2\right]dt
\end{equation}
subject to
\begin{gather*}
{_0^CD}_t^{\alpha(t)}x(t)=e^{x(t)}+2e^tu(t),
\quad 0<\alpha(t)\leq 1,\\
x(0)=0.
\end{gather*}
The exact optimal state and control functions are given by
\begin{equation*}
x(t)=t^2,\quad u(t)=\frac{1}{\Gamma(3-\alpha(t))}
t^{2-\alpha(t)}e^{-t}-\frac{1}{2}e^{t^2-t},
\end{equation*}
which minimize the performance index $J$ with the minimum value $J=0$.
In \cite{MR3864326}, a numerical method based on the Legendre wavelet 
has been used to solve this problem with $\alpha(t)=1$. For solving 
this problem with $\alpha(t)=1$, according to our methods,
we have $n=1$. In this case, both approaches introduced 
in Section~\ref{sec:3} give the same result. 
By setting $M=1$, we suppose that
\begin{equation*}
x{'}(t)=A^TB(t)=a_1 \left(t-\frac{1}{2}\right)+a_0,
\end{equation*}
where
\begin{equation*}
A=[a_0,a_1]^T \text{ and }B(t)=[1,t-\frac{1}{2}]^T.
\end{equation*}
The operational matrix of variable-order fractional integration is given by
\begin{equation*}
P_t^{1}=\left[
\begin{array}{cc}
t & 0 \\
-\frac{t}{4} & \frac{t}{2} \\
\end{array}
\right].
\end{equation*}
Therefore, we have
\begin{equation}
\label{6.2}
x(t)=A^TP_t^{1}B(t)=a_0 t+\frac{1}{2} a_1 (t-1) t.
\end{equation}
Moreover, using the control-affine dynamical system, we get
\begin{equation}
\label{6.3}
u(t)=\frac{1}{2}e^{-t}\left(A^TB(t)-e^{A^TP_t^{1}B(t)}\right)
=\frac{1}{2} e^{-t} \left(a_1 \left(t-\frac{1}{2}\right)
+a_0-e^{a_0 t+\frac{1}{2} a_1 (t-1) t}\right).
\end{equation}
By substituting \eqref{6.2} and \eqref{6.3} into \eqref{6.1},
using the Gauss--Legendre quadrature for computing $J$, and,
finally, setting
\begin{equation*}
\frac{\partial J}{\partial a_0}=0,
\quad  \frac{\partial J}{\partial a_1}=0,
\end{equation*}
we obtain a system of two nonlinear algebraic equations in terms of
$a_0$ and $a_1$. By solving this system, we find
\begin{equation*}
a_0=1,\quad a_1=2,
\end{equation*}
which gives the exact solution
\begin{equation*}
x(t)=t^2 \text{ and }
u(t)=te^{-t}-\frac{1}{2}e^{t^2-t}.
\end{equation*}
As it is seen, in the case of $\alpha(t)=1$, our approaches give the exact 
solution with $M=1$ (only two basis functions) compared to the method introduced 
in \cite{MR3864326} based on the use of Legendre wavelets 
with $\hat{m}=6$ (six basis functions).

Since the optimal state function is a polynomial of degree $2$, 
Approach~I gives the exact solution with $M=1$ for every admissible $\alpha(t)$. 
On the other hand, if $\alpha(t)\neq 1$, then ${_0^CD}_t^{\alpha(t)}x(t)\in H^1(0,1)$. 
Therefore, according to the theoretical discussion and the error bound given 
by \eqref{bound:r}, the numerical solution given by Approach~II converges 
to the exact solution, very slowly, that can be confirmed by the results 
reported in Table~\ref{tab:n1} obtained with $\alpha(t)=\sin(t)$ 
and different values of $M$. Furthermore, by considering a different 
$\alpha(t)$, and by applying the two proposed approaches with $M=5$, 
the numerical results for the functions $x$ and $u$ are displayed 
in Figures~\ref{fig:1} and \ref{fig:2}. Figure~\ref{fig:1} displays 
the numerical results obtained by Approach~I, while Figure~\ref{fig:2} 
shows the numerical results given by Approach~II. 
For these results, we have used
\begin{equation}
\label{eq:alpha:used}
\alpha_1(t)=1,
\quad\alpha_2(t)=\sin(t),
\quad\alpha_3(t)=\frac{t}{2},
\quad\alpha_4(t)=\frac{t}{3}.
\end{equation}
Moreover, the numerical results for the performance index, obtained by our
two approaches, are shown in Table~\ref{tab:1}. It can be easily seen that,
in this case, Approach~I gives higher accuracy results
than Approach~II. This is caused by the smoothness of the exact optimal state function $x$.
\begin{table}[!ht]
\centering
\caption{(Example~\ref{ex1}.) Numerical results obtained by Approach~II 
for the performance index with different $M$ and $\alpha(t)=\sin(t)$.}
\label{tab:n1}
\begin{tabular}{cccccc} \hline
$M$ & $1$ &$2$ &$3$&$2$ & $5$\\ \hline
$J$&$6.80\times 10^{-3}$&$2.33\times 10^{-3}$&$1.76\times 10^{-3}$
&$1.57\times 10^{-3}$ &$1.56\times 10^{-3}$\\ \hline
\end{tabular}
\end{table}
\begin{figure}[!ht]
\centering
\includegraphics[scale=0.75]{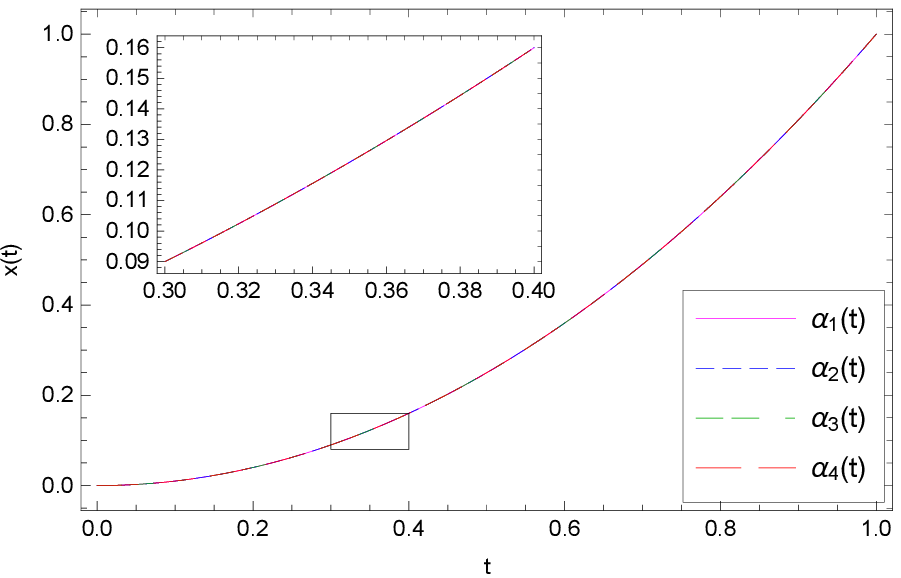}
\quad
\includegraphics[scale=0.75]{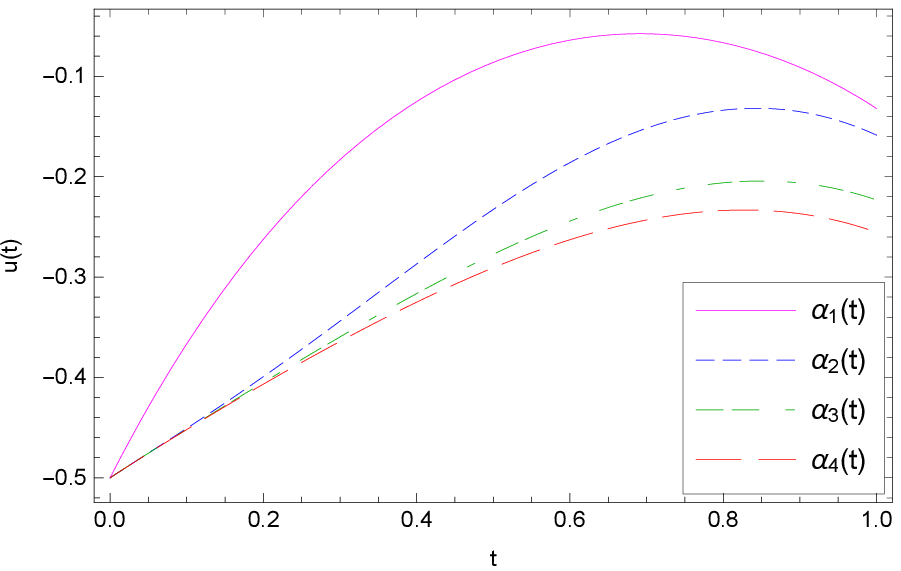}
\caption{(Example~\ref{ex1}.) Comparison between the approximate state (left)
and control (right) functions obtained by Approach~I
with $M=5$ and different $\alpha(t)$ \eqref{eq:alpha:used}.}\label{fig:1}
\end{figure}
\begin{figure}[!ht]
\centering
\includegraphics[scale=0.75]{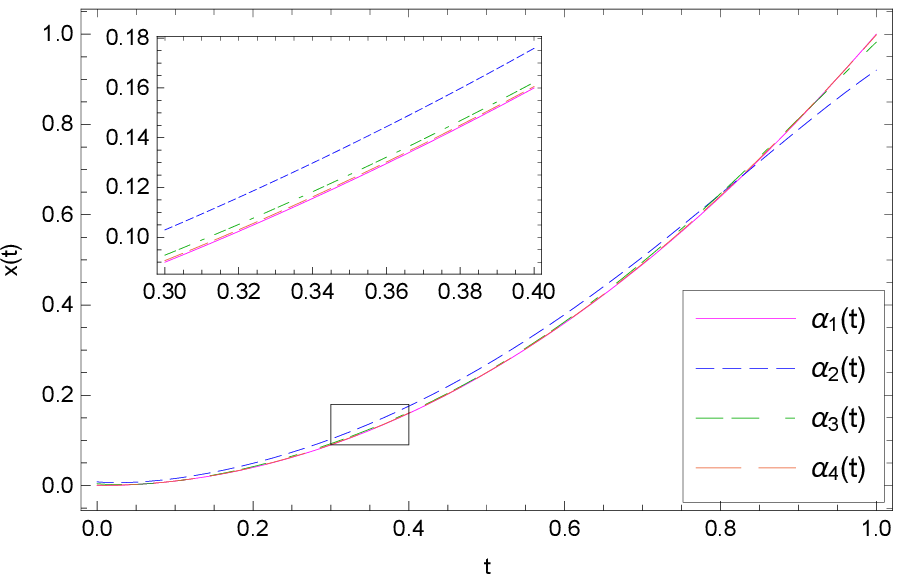}
\includegraphics[scale=0.75]{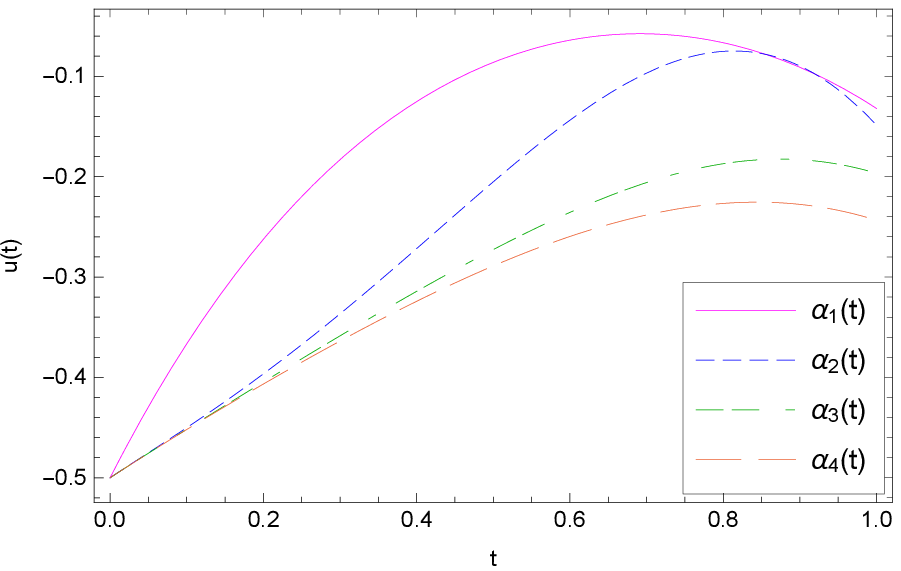}
\caption{(Example~\ref{ex1}.) Comparison between the approximate state (left)
and control (right) functions obtained by Approach~II
with $M=5$ and different $\alpha(t)$ \eqref{eq:alpha:used}.}\label{fig:2}
\end{figure}
\begin{table}[!ht]
\centering
\caption{(Example~\ref{ex1}.) Numerical results for the performance index
with $M=5$ and different $\alpha(t)$ \eqref{eq:alpha:used}.}\label{tab:1}
\begin{tabular}{lcccc} \hline
Method & $\alpha_1(t)$ &$\alpha_2(t)$ &$\alpha_3(t)$&$\alpha_4(t)$\\ \hline
Approach~I
&$3.05\times {10}^{-33}$&$3.26\times {10}^{-33}$
& $6.89\times {10}^{-33}$&$2.08\times {10}^{-33}$\\
Approach~II
&$2.74\times {10}^{-33}$&$1.56\times {10}^{-3}$
&$1.71\times {10}^{-4}$&$2.50\times {10}^{-5}$\\ \hline
\end{tabular}
\end{table}
\end{Example}


\begin{Example}
\label{ex2}
Consider now the following FOC-AP borrowed from \cite{Lotfi}:
\begin{equation}
\label{6.4}
\min~J=\int_0^1\left[\left(x(t)-t^{\frac{5}{2}}\right)^4
+(1+t^2)\left(u(t)+t^6-\frac{15\sqrt{\pi}}{8}t\right)^2\right]dt
\end{equation}
subject to
\begin{equation}
\label{6.5}
{_0^CD}_t^{\frac{3}{2}}x(t)=tx^2(t)+u(t)
\end{equation}
and the initial conditions $x(0)=x'(0)=0$.
For this problem, the state and control functions
\begin{equation*}
x(t)=t^{\frac{5}{2}},
\quad u(t)=-t^6+\frac{15\sqrt{\pi}}{8}t
\end{equation*}
minimize the performance index with the optimal value $J=0$.
We have solved this problem by both approaches. The numerical results
of applying Approach~I to this problem,
with different values of $M$, are presented in Figure~\ref{fig:3}
and Table~\ref{tab:2}. Figure~\ref{fig:3} displays the approximate
state (left) and control (right) functions obtained by $M=1,3,5,7$,
together with the exact ones, while Table~\ref{tab:2} reports the
approximate values of the performance index. Here, we show that
Approach~II gives the exact solution by considering $M=1$.
To do this, we suppose that
\[
{_0^CD}_t^{\frac{3}{2}}x(t)=A^TB(t)=a_1 \left(t-\frac{1}{2}\right)+a_0
\]
with
\begin{equation*}
A=[a_0, a_1]^T \text{ and } B(t)=\left[1, t-\frac{1}{2}\right]^T.
\end{equation*}
Therefore, we have
\begin{equation}
\label{6.6}
x(t)=A^TP_t^{\frac{3}{2}}B(t)=\frac{2}{3 \sqrt{\pi }}(2a_0-a_1)t^{\frac{3}{2}}
+ \frac{8 }{15 \sqrt{\pi }}a_1t^{\frac{5}{2}},
\end{equation}
where
\[
P_t^{\frac{3}{2}}=\left[
\begin{array}{cc}
\frac{4 }{3 \sqrt{\pi }}t^{\frac{3}{2}} & 0 \\
-\frac{2 }{5 \sqrt{\pi }}t^{\frac{3}{2}}
& \frac{8 }{15 \sqrt{\pi }}t^{\frac{3}{2}} \\
\end{array}
\right].
\]
Using the dynamical control-affine system given by \eqref{6.5}, we get
\begin{equation}
\label{6.7}
\begin{split}
u(t)&=a_1 \left(t-\frac{1}{2}\right)+a_0
-t\left(\frac{2}{3 \sqrt{\pi }}(2a_0-a_1)t^{\frac{3}{2}}
+ \frac{8 }{15 \sqrt{\pi }}a_1t^{\frac{5}{2}}\right)^2\\
&=-\frac{64 a_1^2 }{225 \pi }t^6+\left(\frac{32 a_1^2 }{45 \pi }
-\frac{64 a_0 a_1 }{45 \pi }\right)t^5+\left(\frac{16 a_0 a_1 }{9 \pi }
-\frac{16 a_0^2 }{9 \pi }-\frac{4 a_1^2}{9 \pi }\right) t^4+a_1 t
+a_0-\frac{a_1}{2}.
\end{split}
\end{equation}
By substituting \eqref{6.6} and \eqref{6.7} into \eqref{6.4}, the value
of the integral can be easily computed. Then, by taking into account
the optimality condition, a system of nonlinear algebraic equations
is obtained. Finally, by solving this system, we obtain
\begin{equation*}
a_0=\frac{15 \sqrt{\pi }}{16},\quad  a_1=\frac{15 \sqrt{\pi }}{8}.
\end{equation*}
By taking into account these values in \eqref{6.6} and \eqref{6.7},
the exact optimal state and control functions are obtained. Lotfi et al. 
have solved this problem using an operational matrix technique based 
on the Legendre orthonormal functions combined with the Gauss quadrature rule. 
In their method, the approximate value of the minimum performance index 
with five basis functions has been reported as $7.82\times 10^{-9}$ 
while our suggested Approach~II gives the exact value only with two basis functions.
\begin{figure}[!ht]
\centering
\includegraphics[scale=0.75]{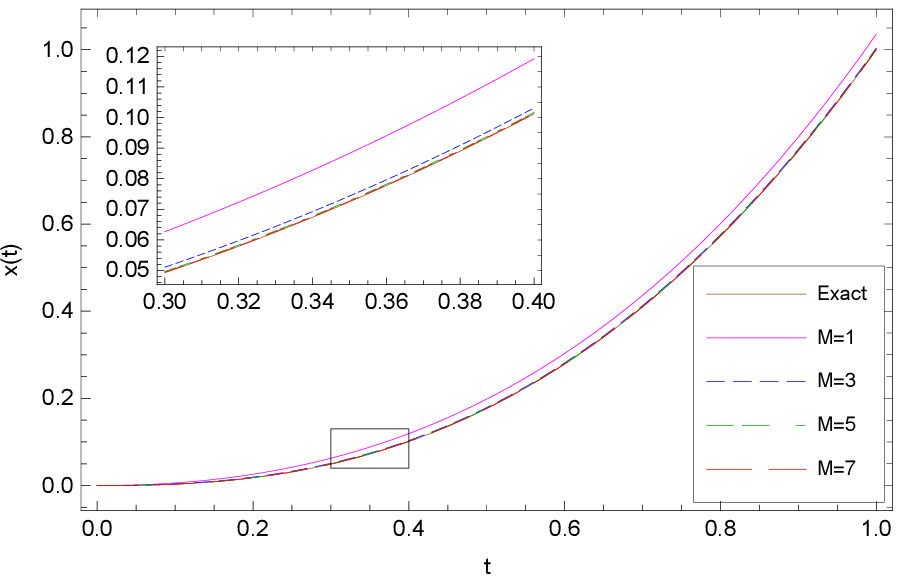}
\quad
\includegraphics[scale=0.75]{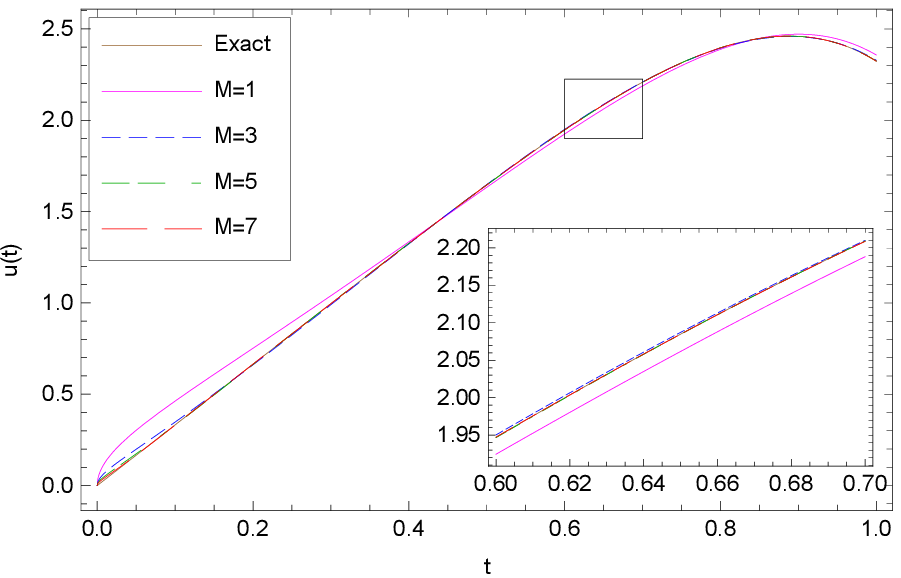}
\caption{(Example \ref{ex2}.) Comparison between the exact and approximate
state (left) and control (right) functions obtained by Approach~I
with different values of $M$.}\label{fig:3}
\end{figure}
\begin{table}[!ht]
\centering
\caption{(Example~\ref{ex2}.) Numerical results for the performance index
obtained by Approach~I with different $M$.}\label{tab:2}
\begin{tabular}{lcccc} \hline
$M$ & $1$ &$3$ &$5$&$7$\\ \hline
$J$&$5.24\times {10}^{-4}$&$7.59\times {10}^{-6}$
& $4.65\times {10}^{-7}$&$5.86\times {10}^{-8}$\\ \hline
\end{tabular}
\end{table}

As we see, in this example, Approach~II yields
the exact solution with a small computational cost, while
the precision of the results of Approach~I
increases by enlarging $M$. Note that here the optimal state
function is not an infinitely smooth function.
\end{Example}


\begin{Example}
\label{ex3}
As our last example, we consider the following FOC-AP \cite{Lotfi}:
\begin{equation*}
\min~J=\int_0^1\left[e^t\left(x(t)-t^{4}+t-1\right)^2
+(1+t^2)\left(u(t)+1-t+t^4-\frac{8000}{77\Gamma\left(\frac{1}{10}\right)}
t^{\frac{21}{10}}\right)^2\right]dt
\end{equation*}
subject to
\begin{eqnarray*}
{_0^CD}_t^{1.9}x(t)=x(t)+u(t),\\
x(0)=1,\quad x'(0)=-1.
\end{eqnarray*}
For this example, the following state and control functions minimize
the performance index $J$ with minimum value $J=0$:
\begin{equation*}
x(t)=t^4-t+1,\quad u(t)=-t^4
+\frac{8000}{77\Gamma\left(\frac{1}{10}\right)}t^{\frac{21}{10}}+t-1.
\end{equation*}
This problem has been solved using the proposed methods with different values of $M$.
By considering $M=1$, the numerical results of Approach~I are shown
in Figure~\ref{fig:4}. In this case, an approximation of the performance index
is obtained as $J=7.21\times {10}^{-1}$. By choosing $M=2$, according to our numerical
method, we have $n=2$. Therefore, we set
\begin{equation*}
x{''}(t)=A^TB(t),
\end{equation*}
where
\begin{equation*}
A=[a_0,a_1,a_2],\quad B(t)=\left[1,t-\frac{1}{2},t^2-t+\frac{1}{6}\right]^T.
\end{equation*}
Hence, using the initial conditions, the state function can be approximated by
\begin{equation*}
x(t)=A^TP^{2}_tB(t)-t+1=\frac{a_2}{12}t^4+\frac{a_1-a_2}{6}t^3
+\left(\frac{a_0}{2}-\frac{a_1}{4}+\frac{a_2}{12}\right)t^2-t+1,
\end{equation*}
where
\[
P^{2}_t=\left[
\begin{array}{ccc}
\frac{t^2}{2} & 0 & 0 \\
-\frac{t^2}{6} & \frac{t^2}{6} & 0 \\
\frac{t^2}{36} & -\frac{t^2}{12} & \frac{t^2}{12} \\
\end{array}
\right].
\]
In the continuation of the method, we find an approximation of the control
function $u$ using the control-affine dynamical system. Then, the method
proceeds until solving the resulting system, which yields
\[
a_0=4,\quad a_1=12,\quad a_2=12.
\]
These values give us the exact solution of the problem. 
This problem has been solved in \cite{Lotfi} with five basis functions 
and the minimum value was obtained as $J=5.42\times 10^{-7}$ while our 
suggested Approach~I gives the exact value with only three basis functions.

In the implementation of Approach~II, we consider
different values of $M$ and report the results in Table~\ref{tab:3}
and Figure~\ref{fig:5}. These results confirm that the numerical solutions
converge to the exact one by increasing the value of $M$. Nevertheless,
we see that since the exact state function $x$ is a smooth function,
it takes much less computational effort to solve this problem
by using Approach~I.
\begin{figure}[!ht]
\centering
\includegraphics[scale=0.75]{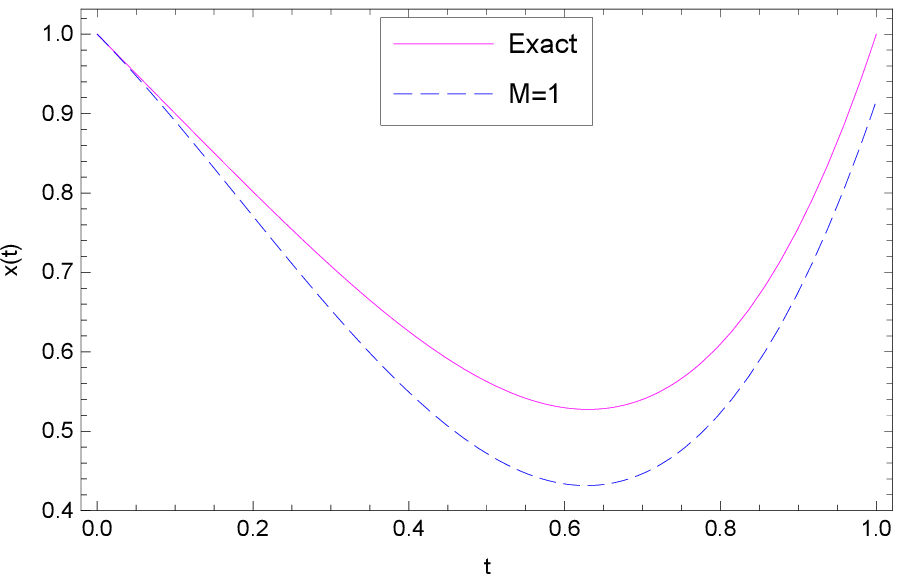}
\quad
\includegraphics[scale=0.75]{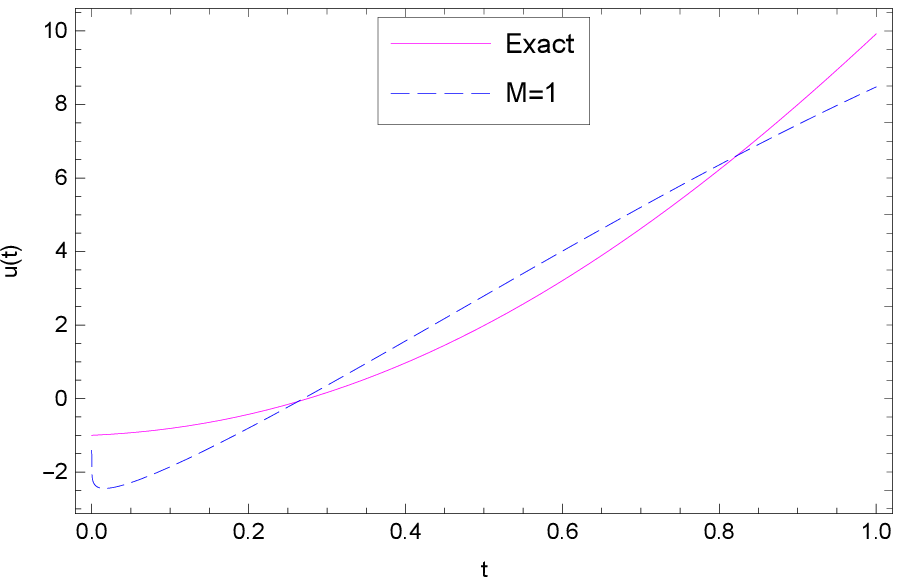}
\caption{(Example~\ref{ex3}.) Comparison between the exact and approximate
state (left) and control (right) functions obtained by Approach~I
with $M=1$.}\label{fig:4}
\end{figure}
\begin{table}[!ht]
\centering
\caption{(Example~\ref{ex3}.) Numerical results for the performance
index obtained by Approach~II with different $M$.}\label{tab:3}
\begin{tabular}{lcccc} \hline
$M$ & $2$ &$4$ &$6$&$8$\\ \hline
$J$&$3.79\times {10}^{-4}$&$5.42\times {10}^{-7}$
& $1.21\times {10}^{-8}$&$7.36\times {10}^{-10}$\\ \hline
\end{tabular}
\end{table}
\begin{figure}[!ht]
\centering
\includegraphics[scale=0.75]{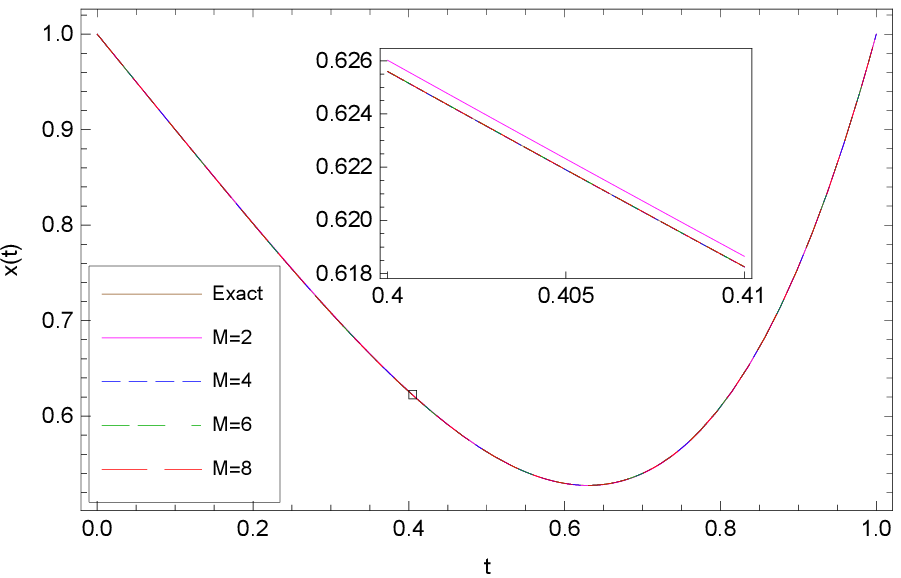}
\quad
\includegraphics[scale=0.75]{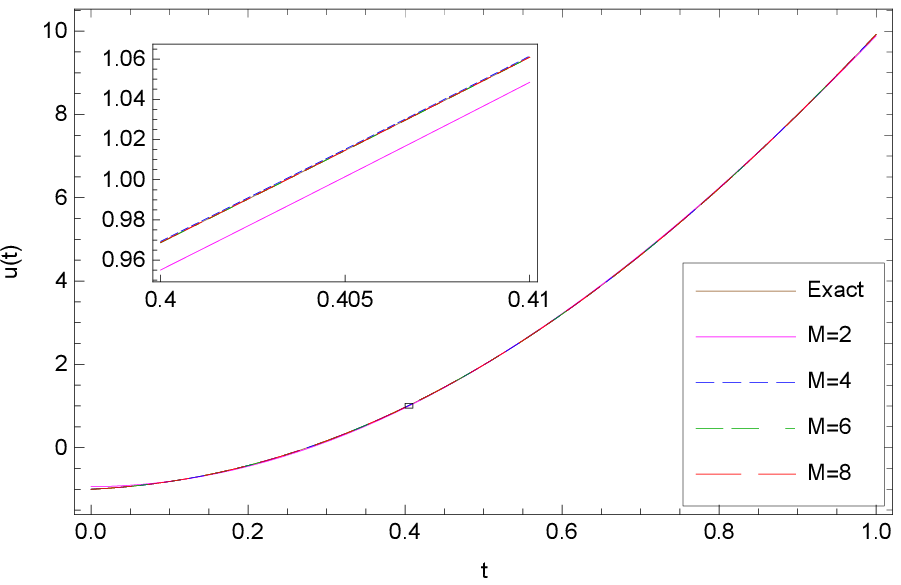}
\caption{(Example~\ref{ex3}.) Comparison between the exact and approximate state (left)
and control (right) functions obtained by Approach~II
with different values of $M$.}\label{fig:5}
\end{figure}
\end{Example}


\section{Conclusions}
\label{sec:6}

Two numerical approaches have been proposed for solving variable-order
fractional optimal control-affine problems. They use an accurate
operational matrix of variable-order fractional integration for Bernoulli polynomials,
to give approximations of the optimal state and control functions. These approximations,
along with the Gauss--Legendre quadrature formula, are used to reduce the original problem
to a system of algebraic equations. An approximation of the optimal performance index
and an error bound were given. Some examples have been solved to illustrate
the accuracy and applicability of the new techniques. From the numerical results
of Examples~\ref{ex1} and \ref{ex3}, it can be seen that our Approach~I
leads to very high accuracy results with a small number of basis functions for optimal control
problems in which the state function that minimizes the performance index is an infinitely
smooth function. Moreover, from the results of Example~\ref{ex2}, we conclude that
Approach~II may give much more accurate results than Approach~I in the cases 
that the smoothness of ${_0^CD}_t^{\alpha(t)}x(t)$ is more than $x^{(n)}(t)$.


\vspace{6pt}

\authorcontributions{Conceptualization, Somayeh Nemati and Delfim F. M. Torres;
Investigation, Somayeh Nemati and Delfim F. M. Torres;
Software, Somayeh Nemati; Validation, Delfim F. M. Torres;
Writing -- original draft, Somayeh Nemati and Delfim F. M. Torres;
Writing -- review \& editing, Somayeh Nemati and Delfim F. M. Torres.}

\funding{Torres was funded by \emph{Funda\c c\~ao para
a Ci\^encia e a Tecnologia} (FCT, the Portuguese Foundation
for Science and Technology) through grant UIDB/04106/2020 (CIDMA).}

\acknowledgments{This research was initiated during a visit of Nemati
to the Department of Mathematics of the University of Aveiro (DMat-UA)
and to the R\&D unit CIDMA, Portugal. The hospitality of the host institution
is here gratefully acknowledged. The authors are also grateful to three
anonymous reviewers for several questions and comments, which helped 
them to improve the manuscript.}

\conflictsofinterest{The authors declare no conflict of interest.}


\reftitle{References}


\end{document}